\documentclass[letterpaper, 10 pt, conference]{ieeeconf}  

\IEEEoverridecommandlockouts                              
\usepackage[para,online,flushleft]{threeparttable}
\usepackage[linesnumbered,ruled,vlined]{algorithm2e}
\usepackage{booktabs}
\usepackage{multirow}
\usepackage{adjustbox}
\usepackage{graphicx}
\usepackage{siunitx}
\usepackage[utf8]{inputenc}
\usepackage[T1]{fontenc}
\usepackage{amsfonts}
\usepackage{gensymb}

\usepackage{amsthm}

\usepackage{epsfig} 
\usepackage{amsmath} 
\usepackage{amssymb}  
\usepackage{siunitx}
\usepackage{subcaption}
\usepackage{caption}
\usepackage{color}
\usepackage[table]{xcolor}
\usepackage{soul}
\usepackage[noadjust]{cite}
\usepackage{amsthm}
\usepackage{mathtools}

\usepackage{hyperref}
\hypersetup{
    colorlinks=true,
    linkcolor=blue,
    filecolor=magenta,      
    urlcolor=blue,
    pdftitle={Overleaf Example},
    pdfpagemode=FullScreen,
    }

\usepackage[acronyms, shortcuts]{glossaries}
\usepackage[noabbrev,capitalise]{cleveref}

\title{\LARGE\bf Task Planning for Multiple Item Insertion using ADMM}

\author{Gavin Zheng$^{1}$
\thanks{$^{1}$Gavin Zheng is with the  Department of Mechanical and Aerospace Engineering, University of California, Los Angeles, CA 90095, USA.
        {\tt\small \{gavinzheng0927\}@ucla.edu}}
}

\begin{document}
\maketitle
\thispagestyle{empty}
\pagestyle{empty}

\begin{abstract}
Mixed-integer nonlinear programmings (MINLPs) are powerful formulation tools for task planning. However, it suffers from long solving time especially for large scale problems. In this work, we first formulate the task planning problem for item stowing into a mixed-integer nonlinear programming problem, then solve it using Alternative Direction Method of Multipliers (ADMM). ADMM separates the complete formulation into a nonlinear programming problem and mixed-integer programming problem, then iterate between them to solve the original problem. We show that our ADMM converges better than non-warm-started nonlinear complementary formulation. Our proposed methods are demonstrated on hardware as a high level planner to insert books into the bookshelf.
\end{abstract}

\section{Introduction} \label{Sec:introduction}
Optimization-based techniques serve as valuable tools in addressing challenges related to robotic motion planning. Various approaches, including mixed-integer convex programs (MICPs) \cite{deits2014footstep, lin2019optimization, dai2019global}, nonlinear or nonconvex programs (NLPs), \cite{DaiValenzuelaEtAl2014,winkler2018gait}, mixed-integer nonlinear programs (MINLPs) \cite{abichandani2015mixed, wang2019trajectory, shin2022design}, and convex optimizations \cite{lin2016convexity, lin2018multi, lin2016formulation} provide robust frameworks for problem formulation. Despite their efficacy, each formulation method has its own set of limitations. Convex optimization, while offering reliable online solutions with guaranteed convergence and rapid solving speed \cite{boyd2004convex}, is constrained to problems amenable to convex formulation. NLPs often converge to local optima, introducing issues such as inconsistent behavior due to reliance on initial guesses. Mixed-integer programs explicitly address discrete variables, commonly utilizing branch-and-bound for solving MIPs \cite{boyd2007branch}. MIP solvers aim for globally optimal solutions, exhibiting more consistent behavior compared to NLP solvers. These solvers typically yield optimal solutions for small-scale problems within a reasonable timeframe \cite{lin2019optimization,tordesillas2019faster}.
. Conversely, MIPs may encounter impractically long solving times when dealing with problems featuring a substantial number of integer variables \cite{lin2021designing}. MINLPs, incorporating both integer variables and nonlinear constraints, offer a high level of expressiveness. Unfortunately, the lack of efficient algorithms for handling MINLPs results in extended solving times and low success rates \cite{lin2022multi}. Consequently, implementing most optimization schemes online becomes challenging.

Recently, researchers have begun exploring the application of machine learning techniques to gather problem-specific heuristics and expedite the optimization solving process. Conventional algorithms employed for solving Mixed-Integer Programs (MIPs), such as branch-and-bound and cutting plane methods, heavily rely on heuristics for efficiently eliminating infeasible regions. Learning methods offer the potential to enhance these heuristics. For instance, in their work, \cite{nair2020solving} utilized graph neural networks to acquire heuristics, while \cite{tang2020reinforcement} employed reinforcement learning to identify effective cutting planes. Alternatively, data collection can be employed to learn and address specific problems or tasks, as demonstrated in studies like \cite{ZhuMartius2020, CauligiCulbertsonEtAl2022, lin2022reduce}. This process effectively transforms into a classification problem, assigning a unique label to each strategy. In the context of online solving, the neural network proposes candidate solutions, effectively reducing the problem to a convex program. Results indicate that CoCo can successfully solve MIPs with approximately 50 integer variables in less than a second. However, it's worth noting that CoCo employs a relatively simple neural network structure. Given that the number of potential strategies grows exponentially with the number of integer variables, the count of unique strategies tends to be close to the total amount of available data.

In numerous multi-agent optimization scenarios, agents may operate in diverse modes defined by discrete variables. The incorporation of discrete mode switching, coupled with collision avoidance constraints, results in a challenging mixed-integer bilinear formulation. To comprehensively assess the efficacy of various data-driven approaches in addressing such problems, \cite{lin2022benchmark} introduced the bookshelf organization problem as a benchmark. This problem involves placing an additional book on a bookshelf with minimal disturbance to existing books. The bookshelf problem serves as an effective benchmark for several reasons: 1) it represents a Mixed-Integer Nonlinear Programming (MINLP) that can be transformed into a Mixed-Integer Convex Programming (MICP) problem with hundreds of integer variables, 2) it can be readily scaled to push algorithms to their limits, and 3) it holds practical significance, allowing for the reasonable collection of data, especially in industries such as logistics.

In this paper, we introduce the Alternating Direction Method of Multipliers (ADMM) to address this problem. ADMM separates the whole MINLP formulation into two parts: a MIP formulation, and a NLP formulation, both solved with applicable solvers. It then iterates between those two formulations until their solutions converge. Despite the process can be slower, it does not require choosing a good warm-start. We demonstrate that this method has a significantly higher chance of finding a solution than the commonly used nonlinear complementary formulation.

\section{Related Works} \label{Sec:related_work}
\subsection{Alternating Direction Method of Multipliers}
The Alternating Direction Method of Multipliers (ADMM) was proposed by \cite{boyd2011distributed} as a general technique to solve optimization problems. The basic idea is to decompose the problem into several subproblems and iterate between them, while using dual variables to help with the convergence. If the magnitude of dual variables are larger, the constraint associated with this dual imposed more constraint force against the objective function, aiding its convergence. ADMM has been introduced into robotics community for real-time control \cite{aydinoglu2022real}, and motion planning \cite{shirai2022simultaneous, lin2022multi}. Despite showing fast convergence on convex problems \cite{stellato2020osqp} and mixed-integer problems \cite{aydinoglu2022real}, its convergence on mixed-integer non-convex problems is slow and not guaranteed.

\subsection{Mixed-integer Programming}
To allow online adaptation, one needs to solve an online optimization formulation in real-time.
This online formulation can have several forms.
Typical robotics problem includes discrete contact and discrete operating modes.
There are generally two ways to include discrete variables in the optimization formulation.
Mixed-integer programs include continuous variables and explicit discrete variables.
The definition of mixed-integer variables is that if the discrete variables are relaxed into continuous ones, the problem becomes convex \cite{dai2019global}.
Typical mixed-integer programming solvers use branch-and-bound methods \cite{boyd2007branch}, cutting plane methods \cite{marchand2003cutting}, and rely on heuristics to choose the branch to explore.
Despite the worst case of solving time, a large portion of problems only requires exploring a small portion of the search tree \cite{williams2013model}.
Mixed-integer programs have been implemented for online motion planning such as \cite{tordesillas2019faster}, and control tasks such as \cite{lin2022learning}.

\subsection{Mathematical Programming with Complementary Constraints}
On the other hand, mathematical programming with complementarity constraints (MPCC) models discrete modes through continuous variables with complementary constraints.
Complementary constraints enforce a pair of variables such that if one of them is non-zero, the other one should be zero.
This constraint is traditionally hard to solve and is sensitive to initial guesses.
Algorithms such as time-stepping \cite{anitescu1997formulating}, pivoting \cite{drumwright2015rapidly}, central path methods \cite{kojima1991unified} are proposed to resolve complementarity.
In the robotics community, complementary constraints are typically used to optimize over gaits for trajectory optimization \cite{PosaCantuEtAl2014, zhang2021transition} or control with implicit contacts \cite{cleac2021fast} where contact forces and distance to the ground are complementary with each other.
On the other hand, complementary constraints can also be used to model binary variables \cite{park2018semidefinite, lin2022developing}. 

\section{Task Planning Problem for Item Insertion} 
\label{Sec:task_planning}
In this section, we outline the task planning problem that involves inserting a sequence of items into a storage bin already occupied by items. The primary goal is to minimize the displacement of the existing items, thereby simplifying the insertion motion executed by the robot manipulator. Consider a 2D shelf with a constrained width $W$ and height $H$ housing rectangular items where each item $i$ possesses dimensions $W_{i}$ (width) and $H_{i}$ (height) for $i=1,...,N$. A series of new items, labeled $i=1,...,N$, is slated for insertion into the shelf in a predetermined order. The shelf accommodates a variety of items in different orientations, potentially requiring the adjustment of existing items to facilitate the insertion of new ones. This study specifically concentrates on the insertion of a single item, utilizing a single-arm robot executing a solitary motion, during which the item remains consistently gripped.

\subsection{Collision free constraint between item shapes}
We solve a sequence of insertion motion under iteration variable $t=1,...,T$. The variables that characterize item $i$ are: position $\textbf{x}_{i}[t] = [x_{i}[t], y_{i}[t]]$ and angle $\theta_{i}[t]$ about its centroid. $\theta_{i}[t]=0$ when a item stands upright. The rotation matrix is: $\textbf{R}_{i}[t] = [cos(\theta_{i}[t]), \ -sin(\theta_{i}[t]); \ sin(\theta_{i}[t]), \ cos(\theta_{i}[t])]$. Let the 4 vertices of item $i$ be $\textbf{v}_{i,k}$, $k=1,2,3,4$. The constraint:

\begin{equation}
    \textbf{x}_{i}[t] + \textbf{R}_{i}[t]\textbf{h}_{i,k} = \textbf{v}_{i,k}[t]
\end{equation}

shows the linear relationship between $\textbf{x}_{i}[t]$ and $\textbf{v}_{i,k}[t]$, where $\textbf{h}_{i,k}$ is the constant offset vector from its centroid to vertices. Constraint: 

\begin{equation}
    \textbf{v}_{i,k}[t] \subset \text{Book Shelf}
\end{equation}

enforces that all vertices of all items stay within the shelf, a linear constraint. Constraint: 

\begin{equation}
    \textbf{R}_{i}^{T}[t]\textbf{R}_{i}[t] = \textbf{I}, \ \ det(\textbf{R}_{i}[t]) = 1
\end{equation}

enforces the orthogonality of the rotation matrix. Constraint:

\begin{equation}
    \textbf{R}_{i}(1, 1)[t] \geq 0
\end{equation}

enforces that the angle $\theta_{i}$ stays within $[-90\degree, 90\degree]$, storing items right side up.

To ensure that the final item positions and orientations do not overlap with each other, separating plane constraints are enforced. The two shapes do not overlap with each other if and only if there exists a separating hyperplane $\textbf{a}^{T}\textbf{x}=b$ in between \cite{boyd2004convex}. We use a set of binary variables $z_{ij}$ to describe if two items $i$ and $j$ are in collision where collision happens if $z_{ij}[t]=1$. For any point $\textbf{p}_{1}[t]$ inside shape 1 then $\textbf{a}^{T}[t]\textbf{p}_{1}[t] \leq b[t] - (1-z_{ij}[t])\delta$, and for any point $\textbf{p}_{2}[t]$ inside shape 2 then $\textbf{a}^{T}[t]\textbf{p}_{2}[t] \geq b[t] + (1-z_{ij}[t])\delta$. The additional amount $\delta$ is used such that if collision does not happen, the two items have a minimal distance. Constraint:

\begin{equation}
    \textbf{a}_{j}^{T}[t]\textbf{a}_{j}[t] = 1
\end{equation}

enforces $\textbf{a}[t]$ to be a normal vector. 

\subsection{Static equilibrium}
Each item should satisfy static equilibrium constraint shown by:

\begin{equation}
\begin{aligned}
& \sum_{j} (\textbf{f}_{ij,i}[t] + \textbf{f}_{ji,i}[t]) + \textbf{f}_{g,i}[t] = \textbf{G} \\
& \sum_{j} \textbf{r}_i[t] \times \textbf{f}_{ij,i}[t] + \sum_{j} \textbf{r}_j[t] \times \textbf{f}_{ji,i}[t] \\
& + \sum \textbf{r}_g[t] \times \textbf{f}_{g,i}[t] = 0 \\
\end{aligned}
\end{equation}

This includes force equilibrium and moment balance which is taken with respect to the item geometric center (assumed to be the center of mass). The external forces include contact forces from any other item, contact forces from the ground or the wall.

Finally, we add additional robustness constraint:

\begin{equation}
\text{If} \ \ \text{$\sum_{j}(z_{ij}[t] + z_{ji}[t]) + z_{gi}[t] = 0$}, \ \text{Then} \ \text{$\theta_{i}[t]=0\degree$}
\end{equation}

This is to exclude one specific case: the item balance at an angle with the center of mass right on top of the contact point and without any other contact point.

\subsection{Objective function}
The objective function the movement of original stored items with respect to its stored pose given by $||\textbf{x}_{i}[t]-\textbf{x}_{i,0}||^2_{W_x} + ||\textbf{R}_{i}[t]-\textbf{R}_{i,0}||^2_{W_\theta}$, and the item poses with respect to their original poses before inserting the new item given by $||\textbf{x}_{i}[t]-\textbf{x}_{i}[t-1]||^2_{W_x} + ||\textbf{R}_{i}[t]-\textbf{R}_{i}[t-1]||^2_{W_\theta}$. The $W_x$ and $W_\theta$ are weights. 

Overall, this is a problem with integer variables, ${z_{ij}[t]}$'s and ${z_{gi}[t]}$'s, and non-convex constraints, hence, an MINLP problem. In practical terms, tackling this problem poses significant challenges in obtaining high-quality solutions. When employing robots for item storage, the permissible solving time is limited to several seconds, and suboptimal solutions result in prolonged realization times. Various potential approaches have been considered to address this issue, including fixing a set of nonlinear variables and solving a Mixed-Integer Convex Program (MICP) \cite{posa2015stability}, transforming nonlinear constraints into piece-wise linear constraints and incorporating them into an MICP \cite{deits2014footstep}, or directly applying MINLP solvers like BONMIN \cite{bonami2007bonmin}. As anticipated, these approaches encounter difficulties in meeting the specified requirements. In the next section, we propose ADMM to solve this formulation at high success rate with decently good cost and solving time.

\section{Solving Task Planning using ADMM} 
\label{Sec:learning_algorithm}
The problem in the previous section can be placed inside a more general frame work. We are interested in solving a type of task planning problem to manipulate multiple items with potential contact between any items. The variables include item position and orientation, contact positions and contact forces. This problem can be formulated into mixed-integer bilinear programs. Assume that a set of problems is given which are parametrized by $S$ drawn from a distribution $D(S)$. For each $S$, we seek a solution to the optimization problem:

\begin{equation}
\begin{aligned}
& \underset{\textbf{x}, \ \textbf{z}}{\text{minimize}} \ f_\textrm{obj}(\textbf{x}, \textbf{z}; \boldsymbol{S}) \\
\text{s. t.} \ \ & f_{i}(\textbf{x}, \textbf{z}; S) \leq 0, \ \ i = 1,...,m_{f} \\
& b_{j}(\textbf{x}, \textbf{z}; S) \leq 0, \ \ j = 1,...,m_{b} \\
& \textbf{z} \in \{0, 1 \}^{n_z}
\label{Eqn:General_formulation}
\end{aligned}
\end{equation}

Where $\textbf{x}$ denotes continuous variables and $z \in \{0, 1\}^{n_z}$ the binary decision variables. Constraints $f_{i}$ are mixed-integer convex, meaning if the binary variables $\textbf{z}$ are relaxed into continuous variables $\textbf{z} \in [0, 1]$, $f_{i}$ becomes convex. Constraints $b_{j}$ are mixed-integer nonlinear, meaning that relaxing the binary variables gives nonlinear constraints. Without loss of generality, $\textbf{x}$ and $\textbf{z}$ are assumed to be involved in each constraint.

Broadly speaking, there are two avenues for transforming a mixed-integer nonlinear program into either a mixed-integer linear program or a nonlinear program. To convert it into a mixed-integer linear program, one must transform non-convex constraints into mixed-integer linear constraints. For instance, trigonometry constraints can be converted into piece-wise linear constraints \cite{deits2014footstep}, or nonlinear constraints can be transformed into mixed-integer envelope constraints \cite{dai2019global}. However, it is worth noting that mixed-integer envelope constraints tend to significantly prolong the solving time. In the study presented in \cite{lin2021designing}, a trajectory planner based on a mixed-integer formulation with envelopes took several hours to solve a trajectory for walking just a few steps upstairs. The alternative approach involves converting binary variables into continuous variables with complementary constraints. This method necessitates the solver to handle a potentially large number of complementary constraints, a computationally challenging task. Directly solving the Nonlinear Program (NLP) with complementary constraints, without an informed initial guess, poses a high risk of infeasibility \cite{zhang2021transition}, unless specialized treatments such as those outlined in \cite{anitescu1997formulating} are employed. However, these methods tend to result in slower solving speeds, as indicated by our benchmark results.

In this paper, we propose to solve problem ~\eqref{Eqn:General_formulation} using ADMM. This formulation separates the original MINLP problem into two parts - a MIP formulation and a NLP formulation. It then alternates between them until their solutions reach a consensus. The ADMM formulation is less restricted by initial guesses, hence the success rate of solving is higher. We describe the implementation details for each formulation. 

The ADMM iterates the following 3 steps recursively:

\begin{align}
\textbf{y}_{M}^{i+1} = arg\min_{\substack{\textbf{y}_{M}}} \ \ \mathcal{L}(\textbf{y}_{M}, \textbf{y}_{N}^{i}, \textbf{w}^{i}) \label{Eqn:Aug_eqn1} \\
\textbf{y}_{N}^{i+1} = arg\min_{\substack{\textbf{y}_{N}}} \ \ \mathcal{L}(\textbf{y}_{M}^{i+1}, \textbf{y}_{N}, \textbf{w}^{i}) \label{Eqn:Aug_eqn2} \\
\textbf{w}^{i+1} = \textbf{w}^{i} + \textbf{y}_{M}^{i+1} - \textbf{y}_{N}^{i+1}  \label{Eqn:Aug_eqn3} 
\end{align}

Where $\textbf{y}=(\textbf{x}, \textbf{z})$ are all variables. $\textbf{w}$ denotes the dual variables. The solving process continues iteratively until the specified consensus condition is achieved such that $||\textbf{y}_{M} - \textbf{y}_{N}|| \leq \Delta$, where $\Delta$ is the specified threshold.

In the first step, we solve ~\eqref{Eqn:Aug_eqn1}, which is converted into a mixed-integer formulation in (\ref{Eqn:MIP_formulation}):  

\begin{equation}
\begin{aligned}
& \underset{\textbf{y}_{M}^{i}}{\text{minimize}} \ f_\textrm{obj}(\textbf{y}_{M}^{i}) \\
& ||\textbf{y}_{M}^{i}-\textbf{y}_{N}^{i}+\textbf{w}^{i}||_{\textbf{G}_{k}}\\
\text{s. t.} \ \ & f_{i}(\textbf{y}_{M}^{i}) \leq 0, \ \ i = 1,...,m_{f} \\
& \textbf{z}_{M}^{i} \in \{0, 1 \}^{dim(\textbf{z})}
\label{Eqn:MIP_formulation}
\end{aligned}
\end{equation}

where the indicator functions in the objective function are converted to explicit constraints again. This problem can be solved using off-the-shelf solvers such as Gurobi.

In the second step, we solve ~\eqref{Eqn:Aug_eqn2}, which is converted into a nonlinear formulation in (\ref{Eqn:NLP_formulation}):
\begin{equation}
\begin{aligned}
& \underset{\textbf{y}_{N}^{i}}{\text{minimize}} \ f_\textrm{obj}(\textbf{y}_{N}^{i}) \\ 
& ||\textbf{y}_{N}^{i}-(\textbf{y}_{M}^{i}+\textbf{w}^{i}))||_{\textbf{G}_{k}}\\
\text{s. t.} \ \ & f_{i}(\textbf{y}_{N}^{i}) \leq 0, \ \ i = 1,...,m_{f} \\
& b_{j}(\textbf{y}_{N}^{i}) \leq 0, \ \ j = 1,...,m_{b} \\
& \textbf{z}_{N}^{i} \in [0, 1]^{dim(\textbf{z})}
\label{Eqn:NLP_formulation}
\end{aligned}
\end{equation}

This problem can be solved using off-the-shelf solvers such as IPOPT. 

Finally, we update the dual variables and weights: 
\begin{equation}
\begin{aligned}
    \textbf{w}^{i+1} = \textbf{w}^{i} + \textbf{y}_{M}^{i} - \textbf{y}_{N}^{i} \\
    \textbf{G}_{k+1} = \gamma\textbf{G}_{k} \\
    \textbf{w}^{i+1} = \textbf{w}^{i} / \gamma
\end{aligned}
\label{eqn:update_weights}
\end{equation}

We show the complete process of ADMM as Algorithm \ref{algorithm}

\begin{algorithm}
\label{algorithm}
\caption{ADMM for item insertion}
\label{algorithm}
\KwIn{$\gamma$, $\textbf{w}_{ini}$, $\textbf{G}_{k,ini}$, $\Delta$}
\textit{intialization} $\textbf{w}=\textbf{w}_{ini}$, $\textbf{G}_{k}=\textbf{G}_{k,ini}$\\
\While{$\Delta\textbf{y}_M \leq \Delta$}{
Solve MIP formulation ~\eqref{Eqn:MIP_formulation} with applicable solvers using current $\textbf{y}_N^i$, $\textbf{w}^i$, $\textbf{G}_{k}$, get $\textbf{y}_M^i$\\
Solve NLP formulation ~\eqref{Eqn:NLP_formulation} with applicable solvers using current $\textbf{y}_M^i$, $\textbf{w}^i$, $\textbf{G}_{k}$, get $\textbf{y}_N^i$\\
Update $\textbf{w}^i$, $\textbf{G}_{k}$ using ~\eqref{eqn:update_weights}}
\Return{$\textbf{y}_M$}
\end{algorithm}

\section{Experiment} 
\label{Sec:experiment_setup}
\subsection{Setup}
We populate shelves with 3, 5, and 7 books, with the goal of inserting an additional book. The problem instances are created in a simulated 2-dimensional environment representing books on a shelf. Initially, four books of randomly determined sizes are placed on the shelf, and one is randomly selected for removal, representing the book slated for insertion. Despite the sequence, the starting configuration with four books represents a valid solution to the problem of placing a book on a shelf with three existing books. This ensures the feasibility of all problem instances. All methodologies are tested on 400 randomly sampled bookshelf problems with the same distribution using the previously mentioned method. The nonlinear ADMM method is implemented on the Mixed-Integer Nonlinear Programming (MINLP) formulation by iteratively transitioning between a mixed-integer formulation and a nonlinear formulation, as detailed in the preceding section. To maintain the non-data-driven nature of this approach, manually crafted heuristics are employed for warm-starting the nonlinear formulation, thereby enhancing the feasibility rate. Considerable effort is invested in tuning the weights $\textbf{G}$, $\gamma$, and appropriately scaling the variables to ensure optimal performance.

After the inserting position is determined, a heuristics-based insertion trajectory is generated and executed by our manipulator. This trajectory first push the item against the item next to the inserting position using position control, pushing that item to the side using force control, then insert the item. The force control is realized using admittance control with a 6-axis force-torque sensor measuring the external wrench. The control law is:

\begin{equation}
    \ddot{\textbf{x}} = M_{d}^{-1}( -D_{d}\dot{\textbf{x}} -K_{d}(\textbf{x}-\textbf{x}_{0}) + K_{f}(\textbf{f}_{meas} - \textbf{F}_{ref}))  
\end{equation}

where $\textbf{x}$, $\dot{\textbf{x}}$ and $\ddot{\textbf{x}}$ are position of the end effector and its first and second derivatives. $\ddot{\textbf{x}}$ is used as control input. $\textbf{f}_{meas}$ is the external wrench measured by force-torque sensor. $\textbf{x}_{0}$ is the reference position. $\textbf{F}_{ref}$ is reference wrench. $M_{d}$, $D_{d}$, $K_{d}$, $K_{f}$ are desired mass, damping, spring coefficient and force sensitivity of the system.

\begin{table*}[]
\caption{ADMM performance for various number of items}
\centering
\begin{tabular}{@{}cccccccc@{}}
\toprule
\# Items                 & Method & Success Rate & Avg. Solve Time & Max. Solve Time & Avg. Objective & Avg. \# of iterations & Solver       \\ \hline
\multirow{2}{*}{4 items} & ADMM   & \textbf{94.6\%}       & 3.7 sec         & 15 sec          & 2352           & 21.1                  & Gurobi+IPOPT \\
                         & MPCC   & 19.4\%       & 1.4 sec         & 3 sec           & 1908           & 14.4                  & IPOPT        \\ \hline
\multirow{2}{*}{6 items} & ADMM   & \textbf{91.1\%}       & 5.5 sec         & 22 sec          & 3357           & 44.8                  & Gurobi+IPOPT \\
                         & MPCC   & 12.5\%       & 1.2 sec         & 7.5 sec         & 3002           & 23.8                  & IPOPT        \\ \hline
\multirow{2}{*}{8 items} & ADMM   & \textbf{83.5\%}       & 8.1 sec         & 41 sec          & 4108           & 62.4                  & Gurobi+IPOPT \\
                         & MPCC   & 5.8\%        & 1.6 sec         & 11.6 sec        & 3598           & 47.6                  & IPOPT        \\ \hline
\end{tabular}
\label{Tab:Results}
\end{table*}

\subsection{Results}
\begin{figure*}[t!]
    \centering
    \includegraphics[height=0.3\textwidth, width=0.9\textwidth]{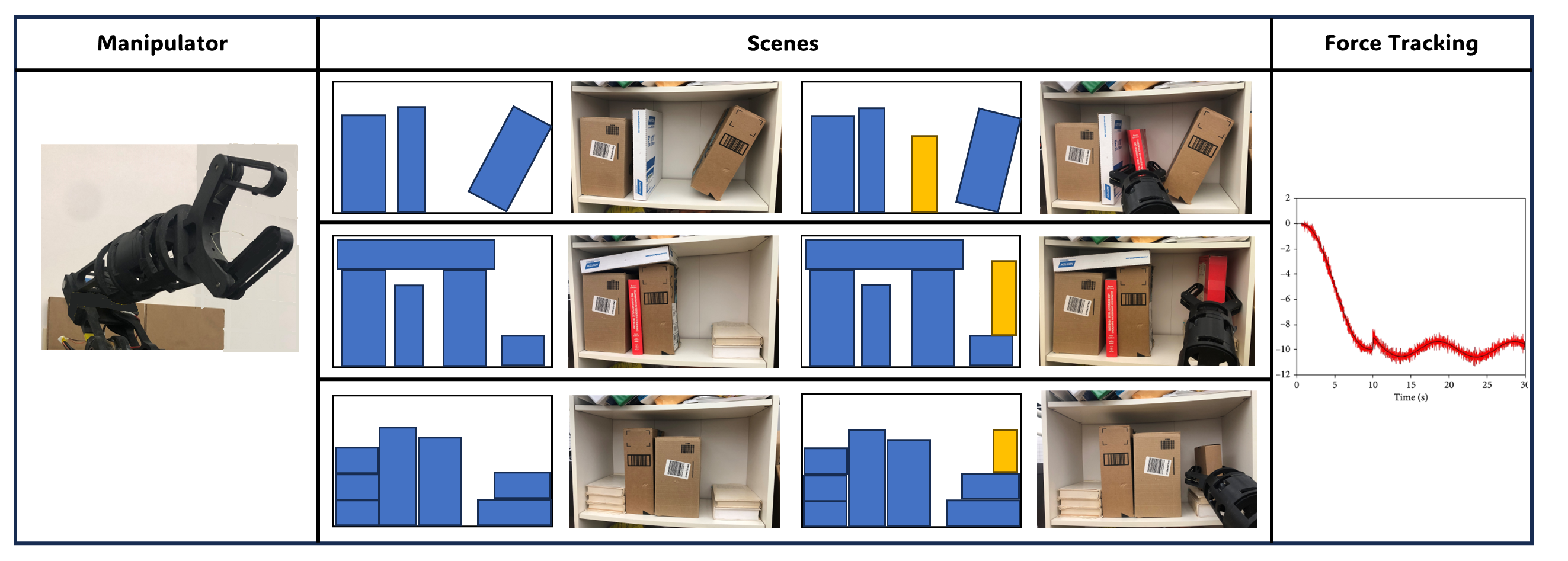}
     \caption{}
     \label{fig:experiment}
\end{figure*}

We implemented Algorithm \ref{algorithm}. The ADMM computation results averaged over the 400 random cases are shown in Table \ref{Tab:Results}.  The outcomes indicate that ADMM performs quite effectively in obtaining feasible solutions. Even without an extensive effort in designing an optimal warm-start, ADMM achieves success rates exceeding 90\% for scenarios involving 4 and 6 items. However, the success rate drops for 8 items, attributed to the increased scale of the problem. In contrast, the MPCC formulation exhibits a success rate of less than 20\%. For nonlinear formulations with complementary constraints, the solver struggles without well-designed warm-starts. The solving speed of ADMM is slower compared to the MPCC formulation. One reason is that ADMM usually takes more than 10 iterations (1 iteration includes 1 MIP and 1 NLP) to get a feasible solution. The nonlinear formulation takes more than 70\% of the solving time, due to the separating plane constraints. If data-driven methods can be used to warm-start the nonlinear formulation, a faster solving speed may be achieved. In addition, the average optimal cost from ADMM is worse. We suspect that the reason is due to the extra consensus terms in the objective function which guides the convergence.

For hardware implementation, we selected 3 scenes and implemented the insertion using the heuristics-based trajectory generator as described in previous section. We used a 6-DoF manipulator equipped with our customized backdrivable actuators. Additional FT sensors are used for high-precision force measurement to validate our proposed controller. This low-level control runs on the hardware at 100Hz. Fig. \ref{fig:experiment} shows the scene before insertion, the scene after insertion, the reference force trajectory along with force measurement from FT sensor. Despite the challenges from contact, friction, and item deformation that are difficult to model and control, this heuristics-driven trajectory achieves the goal in about 30\% of the scenarios we have tested.

\section{Conclusion, Discussion and Future Work} 
\label{Sec:conclusion}
This paper proposes ADMM method to solve a type of mixed-integer nonlinear programmings incorporating item insertion. Numerical experiment and hardware implementation are conducted to support our results. According to the test, ADMM can simultaneously achieve decent optimality, solving speed, and high success rate. Alternative formulations such as MPCC and mixed-integer formulation requires good initial guess to achieve better convergence. Leveraging on pre-solved data, their performances are expected to be similar to our ADMM formulation \cite{er2022data, fletcher2004solving}.

Our future work will be showing that the proposed ADMM formulation can deal with long horizon tasks such as inserting multiple items in sequence. When the problem scales up, current ADMM formulation is expected to perform worse, hence data-driven methods such as \cite{lindeveloping} will be helpful.


{
\bibliographystyle{IEEEtran}
\bibliography{references}
}

\end{document}